# ON THE EXISTENCE OF RIGID $\aleph_1$–FREE ABELIAN GROUPS OF CARDINALITY $\aleph_1$

RÜDIGER GÖBEL AND SAHARON SHELAH

1. INTRODUCTION

An abelian group is said to be $\aleph_1$–free if all its countable subgroups are free. A crucial special case of our main result can be stated immediately.

*Indecomposable $\aleph_1$–free abelian groups of cardinality $\aleph_1$ do exist.*

The first example of any $\aleph_1$–free group which is not free is the Baer–Specker group $\mathbb{Z}^\omega$, which is the cartesian product of countably many copies of the group $\mathbb{Z}$ of integers, known for almost sixty years; cf. Baer [1] or [14, p.94]. Assuming CH, this group of cardinality $2^{\aleph_0} = \aleph_1$ is an example of a non–free abelian group of cardinality $\aleph_1$. Under the same set–theoretic assumption of the continuum hypothesis it can be shown that any countable ring $R$ with free additive group can be realized as the endomorphism ring of an $\aleph_1$–free abelian group $G$ of cardinality $\aleph_1$. The chronologically earlier realization theorem of this kind uses the weak diamond prediction principle which follows from $2^{\aleph_0} < 2^{\aleph_1}$, cf. Devlin and Shelah [6] for the weak diamond, Shelah [28] for the case $\text{End } G = \mathbb{Z}$ and Dugas, Göbel [7] for the case $R = \text{End } G$ and extensions to larger cardinals. Using, what is called Shelah's Black Box, the existence of $\aleph_1$–free groups $G$ with $|G| = 2^{\aleph_0}$ also follows from Corner, Göbel [5] using Dugas, Göbel [7] and combinatorial fine tuning from Shelah [29].

Without the assumption of CH, the existence of non–free, $\aleph_1$–free groups of cardinality $\aleph_1$ follows from a more general result by Griffith [18], Hill [21], Eklof [11], Mekler [24] and Shelah in Eklof [12, p.82,Theorem 8.8]. By an induction it can be shown, that there are $\aleph_n$–free groups, non–free of cardinality $\aleph_n$. The non–abelian case is due to Higman [19, 20].

By Shelah's singular compactness theorem it is known that $\lambda$–free abelian groups of cardinality $\lambda$ do not exist if $\lambda$ is singular, e.g. if $\lambda = \aleph_\omega$, cf. Eklof, Mekler [13]. Hence induction breaks down and it is more complicated to show the existence of $\lambda$–free, non–free abelian groups of cardinality $\lambda > \aleph_\omega$. This question is investigated in Magidor, Shelah [23] and we just refer to this paper and restrict ourselves to cardinals $\lambda \le 2^{\aleph_0}$ again, and we will focus on $\lambda = \aleph_1$. Only very little is known about algebraic properties of $\aleph_1$–free groups of cardinality $\aleph_1$, see Eklof [11] and Eklof, Mekler [13]. Shelah's construction [27] (see also [30]) of groups also mentioned in [12, 13] which are not separable was refined in Eda [10] prove the existence of an $\aleph_1$–free group $G$ of cardinality $\aleph_1$ such that $\text{Hom}(G, \mathbb{Z}) = 0$, a result derived

Part of the work for this paper was carried out while the first author visited Rutgers University. He would like to thank the organizers of MAMLS and the Department of Mathematics for their support.

Number of publication 519. Research was supported by the Edmund Landau Center for research in Mathematical Analysis, supported by the Minerva Foundation (Germany).

independently but later by Corner, Göbel [5]. Counterexamples for Kaplansky's test problems among $\aleph_1$–free groups of cardinality $\aleph_1$ are given recently in Göbel, Goldsmith [17], realizing rings modulo some large ideal, see also [16]. Moreover, $\aleph_1$–separable groups of cardinality $\aleph_1$ serving as counterexamples of Kaplansky's test problems were constructed in [31]. These results about $\aleph_1$–free groups become special cases of our quite satisfying main theorem.

**Main Theorem 4.1** *If $R$ is a ring with $R^+$ free and $|R| < \lambda \leq 2^{\aleph_0}$, then there exists an $\aleph_1$–free abelian group $G$ of cardinality $\lambda$ with* End $G = R$.

We have identified $R$ with endomorphisms acting on the $R$–module $G$ by scalar multiplication. This result has many applications. If $R = \mathbb{Z}$, we derive the existence of $\aleph_1$–free abelian groups of cardinality $\aleph_1$, a result which was unknown.

If $\Gamma$ is any abelian semigroup, then we use Corner's ring $R_\Gamma$, implicitly discussed in Corner, Göbel [4], and constructed for particular $\Gamma's$ in [3] with special idempotents (expressed below), with free additive group and $|R_\Gamma| = \max\{|\Gamma|, \aleph_0\}$. If $|\Gamma| < 2^{\aleph_0}$, we may apply the main theorem and find a family of $\aleph_1$–free abelian groups $G_\alpha(\alpha \in \Gamma)$ of cardinality $\aleph_1$ such that for $\alpha, \beta \in \Gamma$,

$$G_\alpha \oplus G_\beta \cong G_{\alpha+\beta} \quad \text{and} \quad G_\alpha \cong G_\beta \text{ if and only if} \quad \alpha = \beta.$$

Observe that this induces all kinds of counterexamples to Kaplansky's test problems for suitable $\Gamma's$. If we consider Corner's ring in [2], see Fuchs [15, p.145], then it is easy to see that $R^+$ is free and $|R| = \aleph_0$. The particular idempotents in $R$ and our main theorem provide the existence of an $\aleph_1$–free superdecomposable group of cardinality $\aleph_1$, which was unknown as well. Recall that a group is superdecomposable if any non–trivial summand decomposes into a proper direct sum.

Finally, we remark that as the reader might suspect, it is easy to replace $G$ in Theorem 4.1 by a rigid family of $2^\lambda$ such groups with only the trivial homomorphism between distinct members. The main theorem cannot be generalized, replacing $\aleph_1$ by another cardinal. In Section 5 we will show that there are many models of ZFC (e.g. assuming MA and $\aleph_2 < 2^{\aleph_0}$) in which no $\aleph_2$–free group of cardinality $< 2^{\aleph_0}$ has endomorphism ring $\mathbb{Z}$; it is even possible that all such groups are separable and the best one can do now is a realization theorem of the form End $G = R \oplus \text{Ines } G$ with Ines $G \neq 0$ an ideal containing all endomorphisms of finite rank.

This is in contrast to the result [7], that under $\diamond_\lambda$ any countable ring $R$ with $R^+$ free is of the form $R \cong$ End $G$ for all uncountable regular, not weakly compact cardinal $\lambda = |G| > |R|$ such that $G$ is $\lambda$–free. In particular, *the existence of indecomposable $\aleph_2$–free groups of cardinality $\aleph_2$ or the existence of such groups with endomorphism ring $\mathbb{Z}$ is undecidable.*

2. The building blocks, $\aleph_1$–free modules with a distinguished cyclic submodule

Let $R$ be a ring of cardinality $|R| < 2^{\aleph_0}$ such that $R^+$ is a free abelian group. In view of Pontrjagin's theorem we say that an $R$–module is $\aleph_1$–free if any subgroup of finite rank is contained in a free $R$–submodule.

We have the immediate application of Pontrjagin's theorem [14, p.93, Theorem 19.1.].

**Observation 2.1.** *If $M$ is $\aleph_1$–free as $R$–module and $R^+$ is free, then $M$ is $\aleph_1$–free as abelian group, this means all countable subgroups are free.*

**Remark 2.2.** *If $U$ is a finitely generated submodule of an $\aleph_1$–free $R$–module $M$ of infinite rank and $M/U$ is flat, then $M/U$ is an $\aleph_1$–free $R$–module as well.*

**Proof** If $S/U$ is a subgroup of finite rank in $M/U$, then $S_*/U$ denotes its purification and $S_*$ is a pure subgroup of finite rank in $M$, hence it is contained in a free $R$–submodule $F$ of M. Moreover, we find a finitely generated summand $F'$ of the $R$–module $F$ with $S_* \subseteq F'$ and $F/F'$ is $R$–free. Also $F'/U$ is flat because $M/U$ is flat and $F'/U$ can be finitely presented by

$$F'' \to F' \to F'/U \to 0$$

for some finitely generated free module $F''$ mapping onto $U \subseteq F'$. Hence $F'/U$ is projective by Rotman [25, p. 90, 91], and $F/U \cong F/F' \oplus F'/U$ is projective.
Finally we may assume that $F/U$ has infinite rank and $F/U$ is free by a well–known argument of Kaplansky's, cf. [17], for instance. Hence $M/U$ must be an $\aleph_1$–free $R$–module.

Recall that Remark 2.2 does not hold if $U$ is not finitely generated. Consider a free resolution of any torsion–free abelian group A which is not $\aleph_1$–free: $0 \to U \to M \to A \to 0$. By Remark 2.2 in particular quotients of $\aleph_1$–free groups modulo pure, cyclic subgroups are $\aleph_1$–free again.

Next we will construct particular $\aleph_1$–free $R$–modules A with distinguished cyclic submodules $cR$.
First we will fix some more notation. Let $\mathcal{P}$ be a family of $2^{\aleph_0}$ almost disjoint infinite subsets of an infinite set of primes. At present, we choose a fixed $X \in \mathcal{P}$ with an enumeration $X = \{p_n : n \in \omega\}$ without repetitions. Let $T = {}^{\omega>}2$ denote the tree of all finite branches $\eta : n \to 2$, $n < \omega$, where $\ell(\eta) = n$ denotes the length of the branch $\eta$. The branch of length 0 is denoted by $\perp = \emptyset \in T$ and we also write $\eta = (\eta \restriction n - 1)^\wedge \eta(n - 1)$. Finally ${}^\omega 2 = Br(T)$ denotes all infinite branches $\eta : \omega \to 2$ and clearly $\eta \restriction n \in T$ for all $\eta \in Br(T)$, $n \in \omega$.

Let $\lambda$ be an infinite cardinal $\leq 2^{\aleph_0}$ and $Y \subseteq Br(T)$ with $|Y| = \lambda$ and $|R| < \lambda$. Then $V'$ will denote the vector space over the rationals $\mathbb{Q}$ with basis $T \cup Y$. Finally $R$ becomes a vector space by $R \otimes_\mathbb{Z} \mathbb{Q} = \hat{R}$ and $V = V' \otimes_\mathbb{Q} \hat{R}$ is a vector space of dimension $\lambda$. We now select an $R$–submodule $A \subseteq V$ which is generated by $T$ together with elements

$$\eta_0 = \eta,\ \eta_{n+1} = \frac{1}{p_n}(\eta_n + \eta \restriction n + \eta(n) \perp) \in V \qquad (X)$$

defined inductively for all $\eta \in Y$, $n \in \omega$. Hence

$$A = A_X = A_{XY} = \langle \sigma R, \eta_n R : \sigma \in T, \eta \in Y, n \in \omega \rangle \subset V$$

depends on $X \in \mathcal{P}$ and $Y \subseteq Br(T)$. The required cyclic $R$–submodule is $\perp R$. We will show that $(A, \perp R)$ belongs to the category of modules we are interested in, i.e. the following Lemma holds.

**Lemma 2.3.** *Let $(A, \perp R)$ be the pair of $R$–modules defined above, let $B = \langle T \rangle$ and $^- : A \to A/B$ be the canonical homomorphism. Then we have*
*(a) $B$ is a free $R$–module and $A/B = \bigoplus_{\eta \in Y} \bar{\eta}(\bar{X} \otimes_\mathbb{Z} R)$ with $\bar{X} \subseteq \mathbb{Q}$ of characteristic $\chi : \omega \to 2$ with support $X$.*

(b) $A$ is an $\aleph_1$–free $R$–module.
(c) $A/\perp R$ is an $\aleph_1$–free $R$–module.

**Proof** (a) Clearly $B = \bigoplus_{\sigma \in T} \sigma R$ and if $g \in A$, then we use $(X)$ to find $k \in \omega$ and finite sets $T_1 \subseteq T$, $Y_1 \subseteq Y$ with
$$g = \sum_{\eta \in Y_1} \eta_k g_\eta + \sum_{\sigma \in T_1} \sigma g_\sigma$$
for some $g_\eta, g_\sigma \in R$. Using $(X)$ again, we have
$$g \equiv \sum_{Y_1} \eta \frac{g_\eta}{q_k} \mod B$$
where $q_k = \prod_{i=1}^{k-1} p_i$ by the enumeration in $X$ and $\frac{g_\eta}{q_k} \in \bar{X} \otimes R$. Clearly $\{\bar{\eta} : \eta \in Y\}$ is $\mathbb{Q} \otimes R$–independent and hence $\bar{X} \otimes R$–independent and (a) follows.
(b) Obviously $|A| = |Y| = \lambda$. Next we show that
(*) any finite subset of $A$ lies in a submodule $U$ which is free and pure in A.
For any finite subset $E$ of $A$ we can find some $n \in \omega$ and a finite subset $Y_0 \subseteq Y$ such that
$$E \subseteq U = \langle \sigma R, \eta_n R : \sigma \in T, \ell(\sigma) < n, \eta \in Y_0 \rangle.$$
Obviously $U$ is freely generated by the elements $\sigma, \eta_n$. In order to show that $U$ is pure in A, consider $g \in A$ and $m \in \mathbb{N}$ minimal with $gm = u \in U$.
We may write
$$g = \sum_{\eta \in Y_1} \eta_k g_\eta + \sum_{\sigma \in T_1} \sigma g_\sigma \text{ and } u = \sum_{\eta \in Y_2} \eta_n u_\eta + \sum_{\sigma \in T_2} \sigma u_\sigma$$
with $g_\eta, g_\sigma, u_\eta, u_\sigma \in R$ and $k = k(\eta)$ minimal for each $\eta \in Y_1$. Since $gm = u$ we have $Y_1 = Y_2$ and $\eta_n u_\eta = \eta_k g_\eta m$ for all $\eta \in Y_1$ from (a). If $k < n$ for some $\eta \in Y_1 = Y_2$, then we can reduce $Y_1$ to a smaller set $Y_1 \setminus \{\eta\}$ by the observation $\eta_k g_\eta \in U$ and $\eta_k g_\eta m = \eta_n u_\eta$ and $g \in U$ follows by induction. We derive $k \geq n$ for all $\eta \in Y_1$, and suppose $k > n$ for some $\eta$.
We have $p_{k-1}|q = \prod_{i=n}^{k-1} p_i$ and minimality of $m$ requires $p_{k-1}$ does not divide $m$. On the other hand $g_\eta m = q u_\eta$ and $p_{k-1}|q$ hence $p_{k-1}|g_\eta$ which contradicts minimality of $k = k(\eta)$. We derive $k = n$ for all $\eta$ and $g$ decomposes into a $Y$–part $g_Y \in U$ with $g_Y m = \sum_{Y_2} \eta_n u_\eta$ and a $T$–part $g_T \in B$ with $g_T m = \sum_{T_1} \sigma g_\sigma$. However $g_T \in U$, hence $g = g_Y + g_T \in U$ as well and $U$ is pure in A, i.e. (*) holds. Finally $A$ is an $\aleph_1$–free $R$–module by the argument in Remark 2.2 and Pontrjagin's collection of a direct sum of projective modules, see Fuchs [14, p.93, Theorem 19.1.]. Now (b) and also (c) follow from (*).

**Observation 2.4.** *If $(A, \perp R)$ is as above, then $A$ and $A/\perp R$ are $\aleph_1$–free abelian groups with $R \subseteq \text{End } A$, $R \subseteq \text{End }(A/\perp R)$ identifying $r = r \cdot id$ for all $r \in R$.*

Observation 2.4 is immediate from Observation 2.1 and Lemma 2.3, which is all we need in Section 3.
Moreover we will require enough splitting in A which is established by the following

**Proposition 2.5.** *Let $(A, \perp R)$ be as above, where $A = A_X$, $X \neq P \in \mathcal{P}$ and $\bar{P} = \mathbb{Z}_P$ the obvious localization at P. Then $A_X \otimes R_P$ is a free $R_P$–module with $\perp$ a basis element, where $R_P = \mathbb{Z}_P \otimes_{\mathbb{Z}} R$ is the localization of $R$ at $P$.*

**Proof** Recall that
$$A_X = \langle \sigma R, \ \eta_n R : \sigma \in T, \ \eta \in Y, \ n \in \omega \rangle.$$
Moreover $X \cap P$ is finite by our choice of $\mathcal{P}$. We find $k \in \omega$ such that $\{p_n \in X : n \geq k\} \cap P = \emptyset$. Now we claim that
$$T \cup \{\eta_k : \eta \in Y\}$$
is a basis of the $R_P$–module $A_X \otimes R_P$. Note that $\bot \in T$ and Proposition 2.5 will follow.

The set $M = T \cup \{\eta_k : \eta \in Y\}$ is clearly independent over $\mathbb{Q} \otimes_{\mathbb{Z}} R$ in $V$ and hence freely generates the $R_P$–submodule
$$U = \bigoplus_{m \in M} mR_P = F \otimes R_P \subseteq A_X \otimes R_P$$
with $F = \bigoplus_{m \in M} mR$. It remains to show $U = A_X \otimes R_P$.

The submodule $F \subset A_X$ induces a natural sequence
$$0 \to F \to A_X \to A_X/F \to 0$$
of $R$–modules, where $A_X/F$ is generated by $\{\eta_n + F : \eta \in Y, n > k\}$, see Lemma 2.3(a). Using (X) we derive $p_{n-1} \cdot ... \cdot p_{k+1}\eta_n \equiv \eta_k \equiv 0 \mod F$ where the enumeration of primes is taken in X. These primes belong to $\{p_n \in X : n \geq k\}$ and cannot belong to $P$ by our choice of $k$. We observe that $A_X/F$ is a $P'$–group in the well–known sense, that $A_X/F$ is torsion and the order of elements is a product of primes in $P'$, the complement of P. On the other hand $R_P$ is $P'$–divisible, hence $(A_X/F) \otimes R_P = 0$. Using flatness of $R_P$ the above sequence becomes
$$0 \to F \otimes R_P \to A_X \otimes R_P \to (A_X/F) \otimes R_P \to 0$$
and $A_X/F \otimes R_P = 0$ forces $A_X \otimes R_P = F \otimes R_P$ as desired.

## 3. Repeating the building blocks

Let $R$, $\mathcal{P}$ and $|R| < \lambda \leq 2^{\aleph_0}$ be as in Section 2. Then we enumerate $\mathcal{P} = \{X_\alpha : \alpha < \lambda\}$ without repetition and it is easy to find a family $\mathcal{F} = \{L_\alpha \subset \omega : \alpha < \lambda\}$ of infinite, almost disjoint subsets $L_\alpha$ of $\omega$ without repetitions. Since $Br(T) = {}^\omega 2$ and $|{}^\omega 2| = 2^{\aleph_0}$, we can also find a family $\{Y_\alpha \subset Br(T) : \alpha < \lambda\}$ of sets $Y_\alpha$ of branches with the following additional properties

(b1) $|Y_\alpha| = \lambda$ for all $\alpha < \lambda$.
(b2) $Y_\alpha$ has $\lambda$ branch points above every level:
If $\eta \in Y_\alpha$ and $n \in \omega$, there are $\lambda$ distinct branches $\nu \in Y_\alpha$ with $\eta \restriction n = \nu \restriction n$.
(b3) The length of a branch point of branches in $Y_\alpha$ is in $L_\alpha$:
If $\nu \neq \eta \in Y_\alpha$, then $\ell(\nu \cap \eta) \in L_\alpha$.

We use these three families to enumerate a family of $R$–modules $A_{XY}$ constructed in Section 2 defining $A_\alpha = A_{X_\alpha Y_\alpha}$ for all $\alpha < \lambda$. Moreover we denote $R_{X_\alpha} = R_\alpha$ the localization of $R$ at the primes $X_\alpha$ from Section 2.

Inductively we define an ascending, continuous chain of $R$–modules $G_\alpha$ ($\alpha < \lambda$) with distinguished cyclic submodules $c_\alpha R \subset G_\alpha$ for non–limit ordinals $\alpha < \lambda$. The module we are interested in will then be the $R$–module $G = G_\lambda = \bigcup_{\alpha < \lambda} G_\alpha$. If $\alpha = 0$, let $G_0 = \bigoplus_{\nu < \lambda} e_\nu R$ be free $R$–module of rank $\lambda$, which is also a free abelian group of rank $\lambda$ because $R^+$ is free of rank $< \lambda$. We will choose elements $c_\alpha \in G_\alpha$ for non–limit ordinals $\alpha$ subject to the following conditions
(c1) $G_\alpha/c_\alpha R$ is an $\aleph_1$–free $R$–module
(c2) If $c \in G$ and $G/cR$ is an $\aleph_1$–free $R$–module, then $|\{\alpha < \lambda : c = c_\alpha\}| = \lambda$.

The extension $G_{\alpha+1}$ will be constructed such that condition (c1) ensures that $G$ is $\aleph_1$–free and (c2) can easily be arranged by an enumeration of elements $c \in G_\alpha$ with $G_\alpha/cR$ $\aleph_1$–free with $|\alpha|$ repetitions for all $\alpha < \lambda$. If $\alpha = 0$, then for (c1) we may choose a basic element $c_0$ and we do not care for (c2).

If $c_\nu \in G_\nu$ are defined for all $\nu < \alpha$ and $\alpha$ is a limit, then $G_\alpha = \bigcup_{\nu<\alpha} G_\nu$ by continuity and it remains to construct $G_\alpha$ from $c_\beta \in G_\beta$ for $\alpha = \beta + 1$. From our choice (c1) of $c_\beta$ we know that $G_\beta/c_\beta R$ is an $\aleph_1$–free $R$–module. We consider a pushout diagram. There exists a (unique) pushout $R$–module $G_\alpha$ with the well–known pushout mapping properties [14, p.52] or [25] in case of $R$-modules.

$$\begin{array}{ccc} c_\beta R & \longrightarrow & G_\beta \\ \downarrow & & \downarrow \\ A_\beta & \longrightarrow & G_\alpha \end{array}$$

The first row is the canonical embedding and the first column is an embedding by the identification $c_\beta = \bot$. By the pushout property we now may assume that

$(p_\alpha) \qquad G_\alpha = A_\beta + G_\beta$ and $A_\beta \cap G_\beta = c_\beta R$

hence $G_\alpha/c_\beta R \cong G_\beta/c_\beta R \oplus A_\beta/\bot R$. The construction of $G$ is complete.

First we will discuss freeness properties of $G$.

**Lemma 3.1.** *$G$ is an $\aleph_1$–free $R$–module of cardinality $\lambda$.*

**Proof** If $G = \bigcup_{\alpha<\aleph_1} G_\alpha$ as above, then we only have to show that $G_\alpha$ is $\aleph_1$–free for any $\alpha$ which we prove by induction. Since $G_0$ is free we consider $\alpha > 0$ and assume that all $G_\beta$ ($\beta < \alpha$) are $\aleph_1$–free. If $\alpha = \beta + 1$, then $G_\alpha = A_\beta + G_\beta$ and $(p_\alpha)$ holds, hence
$$G_\alpha/c_\beta R \cong G_\beta/c_\beta R \oplus A_\beta/\bot R.$$
The right hand side is $\aleph_1$–free by Lemma 2.3 and assumption on the choice of $c_\beta$. However, if $G_\alpha/c_\beta R$ is $\aleph_1$–free, then $G_\alpha$ must be $\aleph_1$–free as well.

If $\alpha$ is a limit ordinal, then any subgroup of finite rank in $G_\alpha = \bigcup_{\beta<\alpha} G_\beta$ is a subgroup of $G_\beta$ for some $\beta < \alpha$ and $\aleph_1$–freeness follows.

The following observation plays a role in our next proposition, which provides splittings of $G$ coming from Proposition 2.5 and is based on

$R_\alpha \cap R_\beta$ is divisible by all primes not in $X_\alpha \cap X_\beta$ which is finite for $\alpha \neq \beta$.

**Proposition 3.2.** *If $G = \bigcup_{\alpha<\lambda} G_\alpha$ is the $R$–module above, then $G_\alpha \otimes R_\beta$ is a free $R_\beta$–module for all $\alpha \leq \beta < \lambda$.*

**Proof** If $\alpha < \beta$, then $(G_{\alpha+1} \otimes R_\beta)/(G_\alpha \otimes R_\beta) = (G_{\alpha+1}/G_\alpha) \otimes R_\beta$ because $R_\beta$ is a flat $R$–module. We also have $G_{\alpha+1}/G_\alpha = A_\alpha/c_\alpha R$ by the pushout property $(p_{\alpha+1})$ and $(A_\alpha/c_\alpha R) \otimes R_\beta$ is a free $R_\beta$–module by $\alpha \neq \beta$ and Proposition 2.5. We derive that $(G_{\alpha+1} \otimes R_\beta)/(G_\alpha \otimes R_\beta)$ is a free $R_\beta$–module, hence projective and the rest follows inductively by an obvious basis collection. Taking into account that $G_0 \otimes R_\beta$ is a free $R_\beta$–module, the same holds for $G_\alpha \otimes R_\beta$.

**Proposition 3.3.** *With the notation as above we have*
*(a) $A_\beta \otimes R_\beta$ is a direct summand of $G_{\beta+1} \otimes R_\beta$*
*(b) $G_{\beta+1} \otimes R_\beta$ is a direct summand of $G \otimes R_\beta$*
*(c) $A_\beta \otimes R_\beta$ is a direct summand of $G \otimes R_\beta$.*

**Proof** Obviously (c) follows from (a) and (b) and it remains to show the first two assertions.

(a) Observe that $G_\beta \otimes R_\beta$ is free by Proposition 3.2 and we may write $G_\beta \otimes R_\beta = c_\beta R_\beta \oplus H_\beta$ as $R_\beta$–modules by our choice of $c_\beta$. The pushout property $(p_{\beta+1})$ gives $G_{\beta+1} \otimes R_\beta = H_\beta \oplus (A_\beta \otimes R_\beta)$ and (a) follows.

(b) Inductively we will find an ascending, continuous chain of complements $C_\alpha$ of $G_{\beta+1} \otimes R_\beta$ in $G_\alpha \otimes R_\beta$ for $\beta + 1 \leq \alpha \leq \lambda$ and $C_\lambda$ will verify (b). If $\alpha = \beta + 1$, then $C_\alpha = 0$ and if $\alpha$ is a limit ordinal between $\beta + 1$ and $\lambda$ and all $C_\gamma (\gamma < \alpha)$ are defined, then $C_\alpha = \bigcup_{\gamma < \alpha} C_\gamma$ is already defined by continuity and $C_\alpha$ is a complement of $G_{\beta+1} \otimes R_\beta$ in $G_\alpha \otimes R_\beta$ indeed, because $G_\gamma \otimes R_\beta$ ($\gamma \leq \alpha$) is continuous at $\alpha$ as well. It remains to define $C_\alpha$ for $\alpha = \gamma + 1$ where $C_\gamma$ is given. We are in the case $\alpha > \beta + 1$, hence $\gamma > \beta$ and $\gamma \neq \beta$ follows. From Proposition 2.5 we see that $c_\beta R_\beta =\perp R_\beta$ is a summand of the free $R_\beta$–module $A_\gamma \otimes R_\beta$ and we may write $A_\gamma \otimes R_\beta = c_\gamma R_\beta \oplus D_\gamma$. Obviously $C_\alpha = C_\gamma \oplus D_\gamma$ is a complement of $G_{\beta+1} \otimes R_\beta$ in $G_\alpha \otimes R_\beta$ by the pushout property $(p_{\beta+1})$.

## 4. Proof of the Main Theorem

The main result of this paper is the following

**Theorem 4.1.** *If $R$ is a ring with $R^+$ free and $|R| < \lambda \leq 2^{\aleph_0}$, then there exists an $\aleph_1$–free abelian group $G$ of cardinality $\lambda$ with $\mathrm{End}\, G = R$.*

Remark: $G$ will be the $R$–module constructed in Section 3 and we have identified $r \in R$ with $r \cdot id_G$.

**Proof** From Lemma 3.1 we have an $R$–module $G$ of cardinality $\lambda$ which is $\aleph_1$–free as $R$–module, hence $\aleph_1$–free as abelian group. Moreover $R \subseteq \mathrm{End}\, G$ by our identification and we must show that
$\varphi \in \mathrm{End}\, G \setminus R$ does not exist.
Such a homomorphism $\varphi$ has a unique extension $\hat{\varphi} : G \otimes R_\beta \to G \otimes R_\beta$ because $\hat{\varphi} = \varphi \otimes id$ extends and $G \otimes R_\beta / G = (G \otimes R_\beta)/(G \otimes R) \cong R_\beta/R$ being torsion forces uniqueness.
If $c_\alpha \varphi \in c_\alpha R$ for all $\alpha < \lambda$, then $c_\alpha \varphi = c_\alpha r_\alpha$ for some $r_\alpha \in R$. If $\alpha < \lambda$ is fixed, we can choose an element $c \in G$ (even in $G_0$) such that $G/cR$ is an $\aleph_1$–free $R$–module $cR \oplus c_\alpha R$ is a direct sum and $G/(c + c_\alpha)R$ is an $\aleph_1$–free $R$–module as well. There exist some $\gamma, \delta < \lambda$ with $c = c_\gamma$ and $c + c_\alpha = c_\delta$. We have
$c_\gamma r_\gamma + c_\alpha r_\alpha = c_\gamma \varphi + c_\alpha \varphi = (c_\gamma + c_\alpha)\varphi = c_\delta \varphi = c_\delta r_\delta = c_\gamma r_\delta + c_\alpha r_\delta$
and $r_\gamma = r_\delta = r_\alpha$ follows. We find a uniform $r \in R$ such that $c_\alpha \varphi = c_\alpha r$ for all $\alpha < \lambda$. However, $G$ is generated by the set $\{c_\alpha : \alpha < \lambda\}$, hence $\varphi = r$ which was excluded.
There exists $\alpha < \lambda$ such that $c_\alpha \varphi \notin c_\alpha R$. We also find $\gamma > \alpha$ such that $c_\alpha \varphi \in G_\gamma$ and the repetition (c2) (Section 2) of the enumeration of $c_\alpha$'s provides $\gamma < \beta < \lambda$ such that $c_\beta = c_\alpha$, hence
(i) $c_\beta \varphi \notin c_\beta R$ and $c_\beta \varphi \in G_\beta$.
However, $G_\beta \otimes R_\beta$ is a free $R_\beta$–module by Proposition 3.2 and $c_\beta$ is a basic element of the $R_\beta$–module $G_\beta \otimes R_\beta$; we find a free decomposition $G_\beta \otimes R_\beta = c_\beta R_\beta \oplus C$. The pushout $G_{\beta+1} = G_\beta + A_\beta$ gives $G_{\beta+1} \otimes R_\beta = (A_\beta \otimes R_\beta) \oplus C$ and Proposition 3.3(b) provides an $R_\beta$–module $D$ such that $L = C \oplus D$ satisfies

(ii) $(A_\beta \otimes R_\beta) \oplus L = G \otimes R_\beta$, $G_\beta \otimes R_\beta = c_\beta R_\beta \oplus C$
where $C = L \cap (G_\beta \otimes R_\beta)$ by the modular law.
The element $c_\beta \varphi \in G_\beta \subseteq G_\beta \otimes R_\beta$ has a unique decomposition $c_\beta \varphi = c_\beta r + c$ with $r \in R_\beta$ and $c \in C$. If $c = 0$, then $c_\beta \varphi \in c_\beta R_\beta \cap G_\beta = c_\beta R$ by purity of $c_\beta$ is a contradiction. Hence $0 \neq c \in C$ which is a free $R_\beta$–module with a basis B. The element $c = \sum_{b \in [c]} b c_b$ has a unique decomposition and a $B$–support $[c] = \{b \in B : c_b \in R_\beta \backslash \{0\}\} \neq \emptyset$.
On the other hand $c \in C \subseteq G_\beta \otimes R_\beta$ and $cm = \sum_{[c]} bc_b m \in G_\beta \cap C$ for some $m \neq 0$. However $G_\beta \cap C \subset G_\alpha$ for some $\alpha < \beta$, which is contained in the free $R_\alpha$–module $G_\alpha \otimes R_\alpha$. Since $\alpha \neq \beta$, our choice of $R_\alpha, R_\beta$ provides an $h < \omega$ such that
(iii)　　　$p_j$ does not divide $c \in C$ for all $j > h$,
where the enumeration of primes is taken in $X_\beta = \{p_n : n < \omega\}$.
If $\pi : G_{\beta+1} \otimes R_\beta \to C$ denotes the canonical projection induced by (ii), then
(iv)　　　$0 \neq c = c_\beta \varphi \pi$.
Moreover, the image $\eta \varphi \pi$ of any $\eta \in Y_\beta$ viewed as $\eta \in A_\beta \otimes R_\beta \subseteq G_{\beta+1} \otimes R_\beta$ can be expressed by

$$\eta \varphi \pi = \sum_{b \in [\eta]} b r_b^\eta \text{ with } r_b^\eta \in R_\beta \backslash \{0\}$$

with a finite subset $[\eta]$ of $B$. Abusing notation we shall call $[\eta]$ the $B$–support of $\eta$ as well. Recall that $|Y_\beta| = \lambda > |R_\beta| \geq \aleph_0$, and it is easy to find $Y' \subseteq Y_\beta$, $n \in \mathbb{N}$ and $r_b \in R_\beta$ for all $b \in B$ such that $|Y'| = \lambda$ and $|[\eta]| = n$, $r_b^\eta = r_b$ for all $\eta \in Y'$ and $b \in B$. Next we apply the $\Delta$–Lemma to $\{[\eta] : \eta \in Y'\}$ (cf. Jech [22, p.225]) and find $Y'' \subseteq Y'$, $E \subset B$ such that $|Y''| = \lambda$ and $[\eta] \cap [\eta'] = E$ for all $\eta \neq \eta' \in Y''$. Since $[c] \subset B$ is finite, we also find $Y \subset Y''$ such that $|Y| = \lambda$ and $[\eta] \cap [c] \subseteq E$ for all $\eta \in Y$.
From $|Y| = \lambda > \aleph_0$ we find two distinct branches $\eta, \eta' \in Y$ with $\eta \upharpoonright h = \eta' \upharpoonright h$. The branch point $j > h$ of $\eta, \eta'$ belongs to $L_\beta$ by (b3), hence $p_j \in X_\beta$, where $j$ is from the enumeration along branches. The definition branch point gives $\eta \upharpoonright j = \eta' \upharpoonright j$ and $\eta(j) = 1$, $\eta'(j) = 0$ without loss of generality. From the relations $(X_\beta)$ in $A_\beta$ (Section 2) we have $p_j | (\eta'_j + \eta' \upharpoonright j + \eta'(j) \bot)$ in $A_\beta$ and $p_j | (\eta'_j + \eta' \upharpoonright j + \eta'(j) \bot)$ in $A_\beta$, hence $p_j | \eta_j - \eta'_j + \eta(j) \bot = \eta_j - \eta'_j + c_\beta$ in $G_\beta$ and therefore $p_j | (\eta_j \varphi \pi - \eta'_j \varphi \pi) + c_\beta \varphi \pi$. However $[c] = [c_\beta \varphi \pi]$ and if $d = \eta_j \varphi \pi - \eta'_j \varphi \pi$, then $d \upharpoonright E = 0$ by our choice of $\eta, \eta' \in Y$ with $\eta \neq \eta'$, hence $d$ and $c$ are linearly independent. We conclude $p_j | c$ in $C$ which contradicts (iii) and Theorem 4.1 follows.

## 5. A COUNTEREXAMPLE

The reader might suspect that $\aleph_1$ in Theorem 4.1 can be replaced by $\aleph_2$ for instance. This is the case if we assume prediction principles as $\diamondsuit$ (which imply CH), see Dugas, Göbel [7]. However, in general it is no longer true as follows from

**Theorem 5.1.** *Assuming Martin's axiom, any $\aleph_2$–free group of cardinality $< 2^{\aleph_0}$ is separable.*

Recall that an $\aleph_1$–free group is separable if any pure cyclic subgroup is a summand. Preliminaries on (MA) can be seen in Jech [22] or Eklof, Mekler [13].
**Proof**　If $G$ is an $\aleph_2$–free group of cardinality $|G| < 2^{\aleph_0}$, $0 \neq e \in G$ pure in $G$ and $\sigma : e\mathbb{Z} \to \mathbb{Z}$ taking $e\sigma = 1$, then we must extend $\sigma$ to an homomorphism $\Phi : G \to \mathbb{Z}$.

Let $P = \{\varphi; \varphi : D_\varphi \to \mathbb{Z}, e \in D_\varphi, e\varphi = 1\}$ where $D_\varphi$ is a pure and finitely generated subgroup of G. Obviously $|P| < 2^{\aleph_0}$ from $|G| < 2^{\aleph_0}$ and $(P, \subseteq)$ is partially ordered by extensions of maps. Suppose for a moment that $P$ satisfies the hypothesis for MA and $P_g = \{\varphi \in P : g \in D_\varphi\}$ is dense for all $g \in G$. Then by MA there is a compatible set $F \subseteq P$ such that $F \cap P_g \neq \emptyset$ for all $g \in G$. So $\bigcup F = \Phi$ is a partial homomorphism from $G$ to $\mathbb{Z}$. Since $F \cap P_g \neq \emptyset$, also $g \in \mathrm{dom}\ \Phi$ for all $g \in G$, hence $\Phi \in \mathrm{Hom}\ (G, \mathbb{Z})$ and $\Phi$ extends $\sigma$ by definition of P. Thus it remains to show that $(P, \subseteq)$ satisfies the hypothesis of MA:

In order to show that $P_g$ is dense in $P$, we consider any $\varphi \in P$ and find $\varphi \subset \varphi' \in P$ such that $g \in dom\ \varphi'$. Since $G$ is $\aleph_2$–free, there is $D' \supseteq dom\ \varphi$ such that $g \in D'$ and $D'$ is pure and finitely generated in $G$ by Pontrjagin's theorem. Recall that $\mathrm{dom}\ \varphi$ is pure in $G$, hence pure in $D'$ and $D'/\mathrm{dom}\ \varphi$ must be finitely generated and torsion–free. We apply Gauß' theorem to see that $D'/\mathrm{dom}\ \varphi$ is free, hence $D' = \mathrm{dom}\ \varphi \oplus C$ for some $C \subseteq D'$ with $C \cong D'/\mathrm{dom}\ \varphi$. Now it is easy to extend $\varphi$ to a homomorphism $\varphi' : D' \to \mathbb{Z}$. Finally, we must show that $(P, \subseteq)$ satisfies ccc, the countable antichain condition. Let $F \subseteq P$ be an uncountable subset of P. We must find two distinct elements $\varphi_i \in F$ and $\Phi \in P$ such that $\varphi_i \subseteq \Phi$ for $i = 1, 2$. We may assume $|F| = \aleph_1$, hence $(\sum_{\varphi \in F} \mathrm{dom}\ \varphi)_* = U$, the pure subgroup of $G$ generated (purely) by all dom $\varphi$ has cardinality $\aleph_1$ and must be free by hypothesis on G. We select a basis $B$ of $U$ and replace any $\varphi \in F$ by $\varphi'$ with dom $\varphi' = \langle B_\varphi \rangle \supseteq \mathrm{dom}\ \varphi$ with a finite subset of $B_\varphi$ of B. The argument given above allows to extend $\varphi$ to a homomorphism $\varphi'$.

Clearly, it is enough to find two compatible elements $\varphi_i$ in the new F. By the $\Delta$–Lemma (Jech [22, p.225]) we also find $E \subset B$ and $F' \subseteq F$ such that $|F'| = \aleph_1$ and dom $\varphi \cap \mathrm{dom}\ \varphi' = E$ for all $\varphi \neq \varphi' \in F'$. By a pigeon–hole argument we can also find $F'' \subseteq F'$ such that $|F''| = \aleph_1$ and $\varphi \upharpoonright E = \varphi' \upharpoonright E$ for all $\varphi, \varphi' \in F''$. Now it is clear that we can extend two of these maps $\varphi, \varphi'$ to dom $\varphi$ + dom $\varphi'$ as required.

Fachbereich 6, Universität Essen, 45117 Essen, Germany
*E-mail address*: `mat100@vm.hrz.uni-essen.de`

Department of Mathematics, Hebrew University, Jerusalem, Israel, and Rutgers University, Newbrunswick, NJ U.S.A.
*E-mail address*: `shelah@math.huji.ac.il`